\documentclass[12pt,leqno]{amsart}
      \usepackage[english]{babel}
      \usepackage{amsmath}
      \usepackage{amstext}
      \usepackage{amssymb}
      \usepackage{amsthm}
      \usepackage{amsfonts}
      \usepackage[T1]{fontenc}{}{}
      \usepackage{latexsym}
      \usepackage{psfrag}
      \usepackage{setspace}
      
      \theoremstyle{plain}
      \newtheorem{thm}{Theorem}[section]
      \newtheorem*{thm*}{Theorem}
	\newtheorem*{claim*}{Claim}
      \newtheorem{lem}[thm]{Lemma}
      \newtheorem{prop}[thm]{Proposition}

      \theoremstyle{definition}
      \newtheorem{defn}[thm]{\bf Definition}
      
      \newcommand{\N}{\mathbb N}
      \newcommand{\R}{\mathbb R}

      \newcommand{\U}{\mathbb U}

      \newenvironment{appli}{\left( \begin{array}{ccc}}{\end{array} \right)}

      \newcommand{\ba}{\begin{appli}}
      \newcommand{\ea}{\end{appli}}

\begin{document}
\title{Topology of the isometry group of the Urysohn space}
\author{Julien Melleray}
\date{}
\maketitle
\begin{abstract}
Using classical results of infinite-dimensional geometry, we show that the isometry group of the Urysohn space, endowed with its usual Polish group topology, is homeomorphic to the separable Hilbert space $l^2(\N)$.\\
The proof is based on a lemma about extensions of metric spaces by finite metric spaces, which we also use to investigate (answering a question of I. Goldbring) the relationship, when $A,B$ are finite subsets of the Urysohn space, between the group of isometries fixing $A$ pointwise, the group of isometries fixing $B$ pointwise, and the group of isometries fixing $A \cap B$ pointwise.
\end{abstract}

\section{Introduction.}
The Urysohn space $\U$ \footnote{For precise definitions of the objects discussed in this Introduction, see Section 2 below.} was first built by P.S Urysohn more than eighty years ago, see \cite{ury}. It seems to have been considered as little more than a curiosity for a long time, but in the past few decades interest in this space and its isometry group has steadily grown. This is in part thanks to the result, due to V.V. Uspenkij (\cite{Usp2}), that the isometry group of $\U$ is a universal Polish group; it is also the consequence of an  improved understanding of the dynamics and Ramsey-type properties of various ``infinite-dimensional'' structures and their groups of automorphisms, as well as the realization that these structures offer interesting problems and perspectives. At any rate, now the topology of the Urysohn space is completely understood, while its geometry and the properties of its isometry group are the subject of ongoing active research, as witnessed by the volume \cite{Urysohn}.

Uspenskij's result was made possible by Kat\v{e}tov's construction of the Urysohn space, which was based on what we call in this paper \emph{Kat\v{e}tov maps}; these maps provide a good tool to understand one-point metric extensions of a metric space $X$. It is more difficult to manipulate extensions of $X$ by bigger sets; we begin by proving a lemma that shows that one can essentially manipulate metric extensions of a metric space $X$ by sets isometric to a given finite metric space $F$ in much the same way that one manipulates one-point extensions of $X$.

This lemma provides a useful tool for back-and forth constructions where many conditions have to be satisfied all at once; we first use it to prove that any element of the natural basis for the topology of $Iso(\U)$ is homotopically trivial. This in turn implies that $Iso(\U)$ is an ANR, which yields the main result of the paper.

\begin{thm*}
 The isometry group $Iso(\U)$ is homeomorphic to the Hilbert space $\ell_2(\N)$
\end{thm*}
In the last section of the paper, we study stabilizers of finite sets and prove a theorem that illustrates the homogeneity properties of $\U$. This result answers a question that was asked to me by I. Goldbring.\\
Below we make the convention that for $A \subseteq Iso(\U)$, the notation $\langle  A \rangle$ stands for the \emph{closed} subgroup generated by $A$, while if $A\subseteq \U$, $Iso(\U|A)$ denotes the set of isometries fixing $A$ pointwise.

\begin{thm*}
Let $A,B \subseteq \U$ be finite sets. Then
$$Iso(\U |A\cap B)= \langle Iso(\U|A),Iso(\U|B)\rangle  $$
\end{thm*}

For information about Polish spaces and groups, we refer to \cite{Gao} and its bibliography; the volume \cite{Urysohn} is an excellent reference for information on the Urysohn space.
\section{Notations and definitions.}
We recall that a \textit{Polish metric space} is a separable metric space $(X,d)$ whose distance is complete, while a \emph{Polish group} is a separable topological group whose topology admits a compatible complete metric. \\

Whenever $X,Y$ are metric spaces, an \emph{isometry} from $X$ to $Y$ is a distance-preserving map which is also onto. If $f\colon X \to Y$ is distance-preserving but not onto, then we say that $f$ is an \emph{isometric embedding} of $X$ into $Y$. If $(X,d)$ is a Polish metric space we denote by $Iso(X,d)$ (or $Iso(X)$ for short) its isometry group, which we endow with the Polish group topology of pointwise convergence.\\

Let $(X,d)$ be a metric space; a map $f \colon X \to \R$ is said to be a \textit{Kat\v{e}tov map} if it satisfies the following inequalities: 
$$\forall x,y \in X \ |f(x)-f(y)| \le d(x,y) \le f(x)+f(y)\ . $$
These maps identify with one-point metric extensions $X\cup\{z\}$ of $X$, via the correspondence $z \mapsto f=d(z,\cdot)$.\\
We denote by $E(X)$ the set of Kat\v{e}tov maps on $X$; it can be made into a metric space iself by setting, for all $f,g \in E(X)$: 
$$d(f,g) =\sup\{|f(x)-g(x)| \colon x \in X\} \  .$$
The geometric interpretation of that distance is as follows: if $f,g \in E(X)$, then $d(f,g)$ is the smallest possible distance between two points $z_f$, $z_g$ in a metric extension $X\cup\{z_f\} \cup\{z_g\}$ which is such that $d(z_f,x)=f(x)$ and $d(z_g,x)=g(x)$ for all $x\in X$.\\
If $Y \subset X$, then any map $f \in E(Y)$ can be extended to $\hat{f} \in E(X)$, where 
$$\forall x \in X \ \hat{f}(x)=\inf\{f(y)+d(x,y) \colon y \in Y\} \ .$$
The mapping $f \mapsto \hat{f}$ is an isometric embedding of $E(Y)$ into $E(X)$; below we call $\hat{f}$ the \textit{Kat\v{e}tov extension} of $f$.\\
We also need the notion of \textit{support} of a Kat\v{e}tov map $f$: we say that $f \in E(X)$ is \textit{supported} by $S \subset X$ if one has 
$$\forall x \in X \ f(x)=\inf\{f(s)+d(x,s) \colon s \in S\} \ .$$

We recall that the Urysohn space $\U$ is characterized up to isometry, among all Polish metric spaces, by the following property:
$$\forall A \mbox{ finite } \subset \U \ \forall f \in E(A) \ \exists z \in \U \ \forall a \in A \, d(z,a)=f(a)\ . $$
Equivalently, $\U$ is the only Polish metric space which is :
\begin{itemize}
 \item \emph{universal}, i.e any Polish metric space isometrically embeds into $\U$.
 \item \emph{ultrahomogeneous}, i.e any isometry between finite subsets of $\U$ extends to an isometry of the whole space.
\end{itemize}
Below, we will make use of the well-known fact, due to Huhunai\v{s}vili (\cite{huhu}) that $\U$ is \textit{compactly injective}, which means that given any two metric compacta $K \subset L$, every isometric embedding $\varphi \colon K \to \U$ extends to an isometric embedding $\tilde{\varphi} \colon L \to \U$. This implies that an isometry between two compact subsets of $\U$ extends to an isometry of the whole space.\\

Finally, let us make a convention: whenever we manipulate \emph{enumerated} finite metric spaces, we say that $\{a_1,\ldots,b_n\}$ and $\{b_1,\ldots,b_n\}$ are isometric as enumerated finite metric spaces if the map $a_i \mapsto b_i$ is an isometry. If it is clear from the context we will forget to mention the enumeration and just say that $\{a_1,\ldots,a_n\}$ and $\{b_1,\ldots,b_n\}$ are isometric.\\

\section{A lemma about metric spaces.}

As explained in the introduction, the following lemma will be crucial to our constructions. Its aim is to extend the definition of the space $E(X)$ to the case where one considers extensions of $X$ by a given finite metric space; of course the construction yields $E(X)$ as a result if one considers metric extensions of $X$ by a singleton.

\begin{lem}\label{lemlocal1}
Let $X$ be a metric space, and $\{a_1,\ldots,a_n\}$ be an enumerated finite metric space. Then there exists a metric space $Y\supset X$ and subsets $Y_1,\ldots,Y_n \subset Y$ such that $Y= \cup Y_i$, $X=Y_i\cap Y_j$ for all $i \ne j$ and 
\begin{itemize}
 \item[(a)] For each abstract metric extension $X\cup\{a'_1,\ldots,a'_n\}$ of $X$ by an enumerated set isometric to 
$\{a_1,\ldots,a_n\}$, there exist $y_i \in Y_i$ such that $X\cup\{y_1,\ldots,y_n\}$ identifies with $X \cup\{a'_1,\ldots,a'_n\}$ in the natural way.
\item[(b)] $\displaystyle{\forall i \ \forall  y,y' \in Y_i \ d(y,y') = \sup\{|d(y,x)-d(y',x)|\colon x \in X\} \ .}$
\end{itemize}
\end{lem}
One may think of this lemma as a generalization of Kat\v{e}tov's definition of the space $E(X)$: $E(X)$ is the space of extensions of $X$ by one point, the space $Y$ above may be thought of as the set of extensions of $X$ by finite sets isometric to some given finite set. The tricky part is to define a suitable distance, since one must take in consideration $n$-tuples and not points.\\

\noindent {\bf Proof of Lemma \ref{lemlocal1}.} \\
First, pick $n$ disjoint copies of $E(X)$, form their union, and then identify points of $X$ in different copies; this way, one obtains an abstract set $Y=\cup Y_i$ such that each $Y_i$ is naturally identified with $E(X)$, and the intersection of any two distinct $Y_i$'s is $X$.\\
If $f \in E(X)$, we let $f^i$ denote the corresponding element of $Y_i$; we now define a weight function 
$\omega$, by applying the following rules: 
\begin{itemize}
	\item For all $i$, $\omega(f^i,g^i)=d(f,g)$;
	\item For all $i \neq j$, and all $f^i \in Y_i \setminus X$, $g^j \in Y_j \setminus X$, if setting $d(f^i,g^j)=d(a_i,a_j)$ is compatible with 
the triangle inequality for the natural partial distance on $X \cup \{f^i,g^j\}$, then $\omega(f^i,g^j)=d(a_i,a_j)$;
 \item Otherwise, $\omega(f^i,g^j)$ is undefined.
\end{itemize}
To simplify the explanations below, we borrow notation from graph theory and say that a finite sequence $\gamma=(y_0,\ldots,y_n)$ is a \emph{path} joining $y_0$ and $y_n$ if $\omega(y_i,y_{i+1})$ is defined for all $i \leq n-1$. We also say that $y_0,\ldots,y_n$ are the \emph{vertices} of $\gamma$, and define the \emph{weight} $l(\gamma)$ of $\gamma$ by setting
$$l(\gamma)=\sum_{i=0}^{n-1} \omega(y_i,y_{i+1}) \ . $$
We now let $\underline{d}$ denote the pseudometric associated to this weight function; in other words, 
$$\underline{d}(y,y')=\inf\{l(\gamma) \colon \gamma \text{ is a path joining } y \text{ and } y' \} \ .$$
Notice that any two points are joined by a path, so $\underline{d}(y,y')$ is finite for all $(y,y') \in Y^2$; furthermore it is clear that $\underline{d}$ is a distance on $Y$.\\
\begin{claim*}
If $\gamma$ is a path joining $x$ and $y$ and such that $\omega(x,y)$ is defined, then $l(\gamma) \ge \omega(x,y)$; in particular $\underline{d}$ agrees with $\omega$ whenever it's defined.
\end{claim*}
\noindent{\bf Proof of the Claim.} \\
Consider first a path $\gamma$ with three vertices $(x,z,y)$ of elements of $Y$ such that $x \in X$; let us see why 
$\omega(x,z)+\omega(z,y) \geq \omega(x,y).$ This is clear if all three points are in $X$, if $z \in Y \setminus X$ and $y \in X$, and if $z \in X$ and $y \in Y \setminus X$. It is also obvious when $y,z$ both belong to the same $Y_i \setminus X$. So, we may assume that $z=f^i \in Y_i \setminus X$, and $y=g^j \in Y_j \setminus X$, with $i \neq j$. \\
By hypothesis, $\omega(f^i,g^j)$ is well-defined, which implies that $f(x)+d(a_i,a_j) \geq g(x)$ and this is what we needed. \\
Therefore, if $x,x' \in X$, the $\inf$ in the definition of $\underline{d}$ is attained on a path with only two vertices, so that $d(x,x')=\underline{d}(x,x')$. Actually, we have  proved more than that, since we do not need both endpoints of the path to belong to $X$ in order to obtain a path with less vertices. We thus see that $\underline{d}(x,f^i)=f(x)$ for all $x \in X$ and all $f^i \in Y_i$.\\
Now we consider a general path $\gamma=(y_0,\ldots,y_n)$ such that $\omega(y_0,y_n)$ is defined. If $y_0,y_n$ belong to the same $Y_i$ then, by definition of the distance on $E(X)$, $\omega(y_0,y_n)$ is the smallest possible distance $d(z_0,z_n)$ in a two-point metric extension $X\cup\{z_0\} \cup\{z_n\}$ where $d(z_0,x)=d(y_0,x)$ and $d(z_n,x)=d(y_n,x)$ for all $x \in X$.
Since $\underline{d}$ is such a distance, we must have $\omega(y_0,y_n) \le \underline{d}(y_0,y_n) \le l(\gamma)$.\\
We finally consider the case where $y_0=f^i \in Y_i, ,y_n=g^j \in Y_j$, $i\ne j$.\\
If $\gamma$ goes through some $x \in X$, then the reasoning used above shows that $(y_0,x,y_n)$ has to be a path of smaller weight than $l(\gamma)$; the fact that $\omega(y_0,y_n)$ is defined ensures that 
$$\omega(y_0,x)+\omega(x,y_n) \ge d(a_i,a_j)=\omega(y_0,y_n)\ .$$
If $\gamma$ does not go through $X$, then one can first use the preceding case to see that one may assume consecutive vertices of $\gamma$ never belong to the same $Y_i$. Then the triangle inequality in $\{a_1,\ldots,a_n\}$ ensures that $l(\gamma)$ is again more than $d(a_i,a_j)=\omega(y_0,y_n)$.
\hfill $\square_{\text{claim}}$\\

Note that this claim implies that $(Y,\underline{d})$ satisfies condition (b) of the lemma above. Next, we need to see that we indeed embedded in $Y$ all extensions of $X$ by an enumerated set isometric to $\{a_1,\ldots,a_n\}$; for that, it is enough to notice that our construction ensures the following:
$$\forall i \neq j \ \forall f^i \in Y_i \, \forall g^j \in Y_j \ (\omega(f^i,g^j)=d(a_i,a_j)) \Leftrightarrow (\underline{d}(f^i,g^j)=d(a_i,a_j)) $$
We already know that this is true if $f^i $ or $g^j \in X$, so we assume that both elements do not belong to $X$. Then, the implication from right to left is a direct consequence of the definition of $\omega$ and the fact that $\underline{d}$ is a distance. The implication from left to right is an obvious consequence of the fact that $\underline{d}$ extends $\omega$.\\
We have finally done enough to prove (a): if $X \cup\{a'_1,\ldots,a'_n\}$ is an extension of $X$ by an enumerated set isometric to $\{a_1,\ldots,a_n\}$, then each $a'_i$ uniquely defines $y_i \in Y_i$ such that $\underline{d}(y_i,x)=d(a'_i,x)$ for all $i$, and one has $\omega(y_i,y_j)=d(a_i,a_j)$ for all $(i,j)$. Therefore, $\{y_1,\ldots,y_n\}$ is the desired subset of $Y$. \\
This concludes the proof of lemma \ref{lemlocal1}.$\hfill \square_{\text{\ref{lemlocal1}}}$\\

\section{Topology of $Iso(\U)$.}
We prove here that $\mbox{Iso}(\U)$ is homeomorphic to the separable Hilbert space; I am grateful 
to Vladimir Uspenskij for his help and advice on how to prove this result. \\
It should be pointed out here that Uspenskij showed that $\U$ itself is also homeomorphic to $\ell_2(\N)$, see \cite{Usp}. Since $\mbox{Iso}(\U)$ is a Polish group, showing that it is an absolute retract (AR in short) is enough to prove that it is homeomorphic to $\ell_2(\N)$; this is a consequence of a deep result due to Toru\'nczyk and Dobrowolski (\cite{Dobro}).\\ 
And to prove that $\mbox{Iso}(\U)$ is AR, it is enough to show that its topology has a basis that is stable under finite intersections, contains $Iso(\U)$, and is such that for any $V$ in that basis all the homotopy groups $\Pi_n(V)$ are trivial. The reader is invited to consult \cite{VanMill} for proofs and explanations of these facts from infinite-dimensional topology; the results we use here are exposed in chapter 5 of that book.\\

In our study of the topology of $Iso(\U)$, we'll combine Lemma \ref{lemlocal1} with an easy technical lemma, which we state below; if $Y$ is a metric space we let ${\mathcal K}(Y)$ denote the set of compact subsets of $Y$, endowed with the Vietoris topology.

\begin{lem} \label{approx2}
Let $Y$ be a metric space, $X$ be a topological space and $x \mapsto K_x$
be a continuous map from $X$ to ${\mathcal K}(Y)$, $x\mapsto \varphi(x)$ a continuous map from $X$ to $Y$ and $x \mapsto \psi_x$ a continuous map from $X$ to $E(Y)$ such that $\psi_x$ is supported by $K_x$. \\
Then there exists a continuous map $x \mapsto \tau_x$ from $X$ to $E(Y)$ such that $\tau_x(y)=\psi_x(y)$ for all $x \in X$ and $y \in K_x$, and 
$$\tau_x(\varphi(x))= \inf\{f(x) \colon f \in E(Y) \mbox{ and } f(y)=\psi_x(y) \text{ for all } y \in K_x\} \ .$$
\end{lem}

The statement above is somewhat technical, but the idea is simple: given a continuous assignment of partial metric conditions on $Y$, and a continuous assignment of points $\varphi(x)$  in $Y$, one can build a continuous assignment of metric conditions on $Y$ which extends the partial conditions and also takes a minimal value on $\varphi(x)$. \\

\noindent{\bf Proof of lemma \ref{approx2}.} 
Define for all $x \in X$ the map $\tau_x$ by setting successively:
\begin{enumerate} \item $\tau_x(y)=\psi_x(y)$ for all $y\in K_x$;
\item
$\tau_x(\varphi(x))= \sup\{|d(\varphi(x),y)-\psi_x(y)| \colon y \in K_x\}.$
\end{enumerate}
This defines a map $\tau_x$ which belongs to $E(K_x \cup \{\varphi(x)\})$; letting again $\tau_x$ denote its Kat\v{e}tov extension to $Y$, the proof is finished. $\hfill \square_{\text{\ref{approx2}}}$\\

Now we are ready to state the main technical result of this section.

\begin{prop} \label{Pin} Let  $V$ be a basic open set for the topology of $\mbox{Iso}(\U)$, and $K \subseteq L$ be two compact spaces. Let also $\varphi_0$ be some element of $V$, and $\Phi_0,\, \Phi_1 \colon L \to V$ be two continuous maps such that $\Phi_j(k)= \varphi_0$ for $j=0,1$ and all $k \in K$.\\
Then there exists a continuous path $\Phi \colon L \times [0,1] \to V$  between $\Phi_0$ and $\Phi_1$, such that $\Phi(k,t)=\varphi_0$ for all $k \in K$ and all $t \in [0,1]$.
\end{prop}
By "basic open set for the topology of $\mbox{Iso}(\U)$", we mean some element of the natural basis for the product topology.\\
In particular, all the homotopy groups of $V$ are trivial: taking $K=\emptyset$ and $L=\{0\}$ proves that $V$ is path-connected, and letting $K$ be a singleton and $L=\mathbb  S^n$  shows that $\Pi_n(V)$ is trivial for all $n$. \\
As we explained above, this implies that $\mbox{Iso}(\U)$ is AR, hence Proposition \ref{Pin} is enough to prove that $\mbox{Iso}(\U)$ is homeomorphic to the Hilbert space.\\

\noindent{\bf Notation.} If $F,X$ are metric spaces then we denote by $Emb(F,X)$ the set of isometric embeddings 
of $F$ into $X$ (we don't endow it with a topology). \\

Let us now explain the construction in detail; we begin by explaining how to prove proposition 
\ref{Pin} in the case when $V=\mbox{Iso}( \U)$. \\

\noindent{\bf Proof of Proposition \ref{Pin} in the case $V=\mbox{Iso}(\U)$. }\\
In this case, we may assume that $\varphi_0=id$, and that $\Phi_0(l)=id$ for all $l \in L$.\\
We use the back-and-forth method; to apply it we pick a countable dense subset $\{x_n\colon n \in \N\}$ of $\U$.\\
We build an increasing sequence of compact sets $F_n \subset \U$, and maps 
$ (l,t) \mapsto \Phi^n_{l,t}$ from $L \times [0,1]$ to 
$Emb(F_n,\U)$, with the following properties:\\
\begin{enumerate}
\item $\Phi^{n+1}_{l,t}$ extends $\Phi^n_{l,t}$ for all $n$ and all $(l,t) \in L \times [0,1]$; 
\item $\Phi^n_{l,0}(x)=x$ and $\Phi^n_{l,1}(x)=\Phi_1(l)(x) $ for all $x \in F_n$ and all $l\in L$;
\item $\Phi^n_{k,t}(x)=x$ for all $k \in K$ and all $x \in F_n$;
\item $\forall x \in F_n \ (l,t) \mapsto \Phi^n_{l,t}(x)$ is continuous;
\item $\forall n \ x_n \in F_{2n}$; 
\item $\forall n \ \forall l,t  \  \ x_n \in \Phi^{2n+1}_{l,t}(F_{2n+1})$.\\
\end{enumerate}

\noindent  If one lets $F=\cup F_i$, then at the end of this construction we get isometric maps $\Phi_{l,t} \colon F \to \U$, with the property that $(l,t) \mapsto \Phi_{l,t}(x)$ is continuous for all $x$; these maps extend to isometries $\Phi_{l,t} \colon \U \to \U$ with the same property (the construction ensures that they are surjective, and continuity is easy to check). Since we have by construction that $\Phi_{l,0}=\Phi_0(l)$ and $\Phi_{l,1}=\Phi_1(l)$ for all $l \in L$, and $\Phi_{k,t}=\varphi_0$ for all $k \in K$, setting $\Phi(l,t)=\Phi_{l,t}$ defines a map with the desired properties.\\

\noindent We now explain how to carry out the construction above; for this, we need to be able to do two different things:\\
- (forth) If $F$ is a compact set obtained at a previous step of the construction, $(l,t) \mapsto \Phi_{l,t}(x)$ is a continuous map for all $x$ and satisfying all the conditions above, and $z$ is some point not in $F$, then one has to be able to extend the maps $\Phi_{l,t}$ to $F \cup\{z\}$ in such a 
way that they still satisfy all the conditions.\\
- (back) If $F$, $(l,t) \mapsto \Phi_{l,t}$ are as above, and $z\in \U$, then one must find some compact set $F' \supset F$, and a suitable extension of the maps $\Phi_{l,t}$ to $F'$, in such a way that $z \in \Phi_{l,t}(F')$ for all $l,t$.\\
To do this, we first notice that it is easy to find continuously points $z_{l,t}$ such that it would be compatible to set $F'=F \cup\{z_{l,t}\}$, and $\Phi_{l,t}(z_{l,t})=z$; explicitly, the points must satisfy $d(z_{l,t},y)=d(z,\Phi_{l,t}(y))$ for all $l,t$. Then, we have to do a more complicated 
version of the forth step, extending $\Phi_{l,t}$ to each $z_{l,t}$ while preserving continuity. So it should 
be enough to  detail the back step.\\

\noindent {\bf Back step.}\\
Let $ F \subset \U$ be compact, $z \in \U$ and $(l,t) \mapsto \Phi_{l,t}$ be maps satisfying
the conditions for the back step. \\
We begin by picking a continuous map $(l,t) \mapsto z_{l,t}$ such that 
$$\forall y \in F\ d(z_{l,t},y)=d(z,\Phi_{l,t}(y))\ .$$
This is possible because of the compact injectivity of $\U$.\\
Then we fix a countable dense set $\{(l_n,t_n)\}$ in $L\times [0,1]$ and extend inductively the maps $\Phi_{l,t}$ to $F\cup\{z_{l_n,t_n}\}_{n \in \N}$ in such a way that 
\begin{itemize}
\item $\Phi_{l_n,t_n}(z_{l_n,t_n})=z$ for all $n$; 
\item $(l,t) \mapsto \Phi_{l,t}(z_{l_n,t_n})$ is continuous for all $n$;,
\item all the boundary conditions are respected.
\end{itemize} 
Then each $\Phi_{l,t}$ extends 
to $F \cup \{z_{l,t} \colon (l,t) \in L \times [0,1]\}$, and it defines a map $(l,t) \mapsto \Phi_{l,t}$ 
with all the desired properties. \\
We now explain how to carry out this inductive extension; pick $n \in \N$ and assume that we have defined $\Phi_{l,t}(z_{l_0,t_0}),\ldots,\Phi_{l,t}(z_{l_{n-1}t_{n-1}})$ in such a way that, for all $i \leq n-1$ and all $l,t$ one has 
$$d(\Phi_{l,t}(z_{l_i,t_i}),z)=d(z_{l_i,t_i},z_{l,t}) \ .$$ 
Define now
$$ M=F \cup \{z_{l_0,t_0},\ldots,z_{l_{n-1},t_{n-1}}\}, \mbox{ and} $$
$$Y_{l,t}=\{z\}  \cup \{\Phi_{l,t}(m) \colon m \in M\}, \qquad Y= \bigcup_{(l,t) \in L \times [0,1]} Y_{l,t}\ .$$
Let $\psi_{l,t}$ denote the Kat\v{e}tov map on $Y$ with support in $Y_{l,t}$ and values defined by: 
\begin{itemize}
\item $\psi_{l,t} (z) =d(z_{l_n,t_n},z_{l,t})$;
\item $ \psi_{l,t}(\Phi_{l,t}(m))=d(z_{l_n,t_n},m)$ for all $m \in M$.\\
\end{itemize}
Apply now lemma \ref{approx2} to $\psi$ with the functions 
\begin{itemize}
\item $(l,t) \mapsto \varphi_1 (l,t)= z_{l_n,t_n}$;
\item $(l,t) \mapsto \varphi_2(l,t)=z$;
 \item $(l,t) \mapsto \varphi_3 (l,t)  =\Phi_{l,1}(z_{l_n,t_n})$.
\end{itemize}  
This yields three continuous maps $\tau_1,\, \tau_2, \, \tau_3$ with values in $E(Y)$ such that $d(\tau_i(l,t),z)=d(z_{l_n,t_n},z_{l,t})$, 
$d(\tau_i(l,t),\Phi_{l,t}(m))=d(z_{l_n,t_n},m)$ for all $m \in M$, and $\tau_1(l,0)=z_{l_n,t_n}$, $\tau_2(l_n,t_n)=z$ and $\tau_3(l,1)=\Phi_{l,1}(z_{l_n,t_n})$.\\
We are almost done: let $f,g,h \colon L \times [0,1] \to \R$ be three positive-valued maps such that $f+g+h=1$, $f(l,0)=1=h(l,1)$ for all $l \in L$, and $ g(l_n,t_n)=1$. Then, set 
$$z^n_{l,t}=f(l,t)\tau_1(l,t)+g(l,t)\tau_2(l,t)+h(l,t)\tau_3(l,t)_ .$$
Since the map $(l,t) \mapsto z^n_{l,t}$ is continuous, and $\U$ is compactly injective, 
one may assume that $z^n_{l,t}$ belongs to $\U$ for all $(l,t)$. 
Then, setting $\Phi_{l,t}(z_{l_n,t_n})=z^n_{l,t}$ defines a suitable extension of the maps $\Phi_{l,t}$. \hfill $\square_{\text{\ref{Pin}(} V=Iso(\U)\text{)}}$\\

\noindent This concludes the proof of Proposition \ref{Pin} in the case when $V=\mbox{Iso}(\U)$.
The general case is now  not too hard to obtain, thanks to lemma \ref{lemlocal1}: we will use the same 
back-and-forth construction as above, with a special first step that ensures that we stay within a given basic open set $V$.\\

\noindent{\bf Proof of Proposition \ref{Pin} in the general case.}\\
Let $V=\{\varphi \in \mbox{Iso}(\U) \colon d(\varphi(x_i),y_i) < \varepsilon_i (i=1,\ldots n)\}$ be an element of the
 natural basis for the topology of $\mbox{Iso}(\U)$. \\
 Pick some $\varphi_0 \in V$ and continuous maps $\Phi_0,\Phi_1 \colon L \to V$ such that $\Phi_i(k)=\varphi_0$ for all $k \in K$.
This time, we let  
$$Y=\{x_1,\ldots,x_n\}\cup\{y_1,\ldots,y_n\}\cup\{\varphi_0(x_1),\ldots,\varphi_0(x_n)\}\ $$
Take a metric space $Y'$ as in Lemma \ref{lemlocal1} for $Y$, $\{x_1,\ldots,x_n\}$.\\
Then, define for all $(l,t)$ an extension of $Y$ by a set $\{z^{l,t}_É,\ldots,z^{l,t}_n\}$ isometric to $\{x_1,\ldots,x_n\}$ by setting:
\begin{itemize}
\item $d(z_{l,t}^i,z_{l,t}^j)=d(x_i,x_j)$,
\item $\forall y \in Y d(z_{l,t}^i,y)=(1-t)d(\Phi_0(l)(x_i),y)+td(\Phi_1(l)(x_i),y)$ 
\end{itemize} 
One can assume that each $z_{l,t}^i$ belongs to $Y'_i$, and the condition (b) of Lemma \ref{lemlocal1} guarantees that the maps 
$(l,t) \mapsto z_{l,t}^i$ are all continuous, so that $Y \cup\{z_{l,t}^i\}$ is compact.\\
This means as usual that we may assume that $z_{l,t}^i \in \U$, and $d(z^i_{l,t},y_i) <\varepsilon_i$ since it is a convex combination of $d(\Phi_1(l)(x_i),y_i)$ and $d(\Phi_0(l)(x_i),y_i)$, which are both $< \varepsilon_i$. \\
We have just built continuous maps $(t,l) \mapsto z^i_{l,t} \in \U$ 
($i=1,\ldots n$) such that:
\begin{itemize}
\item $z^i_{l,0}=\Phi_0(l)(x_i)$ and $z^i_{l,1}=\Phi_1(l)(x_i)$ for all $l \in L$;
\item $z^i_{k,t}=\varphi_0(x_i)$;
\item $d(z^i_{k,t},z^j_{k,t})=d(x_i,x_j)$;
\item $d(z^i_{l,t},y_i) <\varepsilon_i$ for all $t,l$.
\end{itemize}

One may now set $\Phi_{l,t}(x_i)=z^i_{l,t}$, and proceed with the construction as in the case when $V=\U$ (using the back-and-forth argument that was explained above, starting with the compact set we just built); in the end, we obtain a continuous path $\Phi_{l,t}$ with all the desired properties, the fact that 
$\Phi_{l,t} \in V$ for all $l,t$ being ensured by the beginning of this construction.$\hfill \square_{\text{\ref{Pin}}}$\\

We have finally proved the main result of this section.
\begin{thm}
The group $\mbox{Iso}(\U)$ is homeomorphic to the Hilbert space $\ell_2(\N)$.
\end{thm}

\section{Stabilizers of finite sets.}
In the section we will discuss a result about stabilizers of finite sets in the Urysohn space; the theorem we obtain was conjectured by I. Goldbring.\\

If $A \subset \U$ is a finite set we denote by $Iso(\U|A)$ the set of isometries of $\U$ that coincide with $id_{\U}$ on $A$; if $G \subset Iso(\U)$ we denote by $\langle  G\rangle $ the \emph{closed} subgroup generated by $G$.

\begin{defn}
A metric triangle $\{a,b,c\}$ is called \emph{flat} if one of the triangle inequalities for $\{a,b,c\}$ is actually an equality.
\end{defn}

Below we will use the following observation: assume $X$ is a set, $d_1,\ldots,d_n$ are distances on $X$ and $d=\frac{1}{n}\sum d_i$. Then $d$ is still a distance on $X$, and if a triangle is not flat for some $d_i$ then it is not flat for $d$ either. 

\begin{prop}\label{stab}
Let $A=\{a_1,\ldots,a_n\}$ and $ B=\{b_1,\ldots,b_n\} \subset \U$ be enumerated finite isometric sets, with additionally $a_i=b_i$ for all $i$ less than some $k$ (possibly $0$) and $a_i \ne b_j$ for all $i,j >k$.\\
Let also $\varphi \colon A \to \U $ be a partial isometry such that $\varphi(a_i)=a_i$ for all $i \le k$. \\
Then for any $\varepsilon >0$ there exists $\psi\in \langle Iso(\U|A),Iso(\U|B) \rangle $ such that 
$$\forall i \in\{1,\ldots,n\} \ d(\psi(a_i),\varphi(a_i)) \le \varepsilon \ .$$
\end{prop}

\noindent{\bf Proof of Proposition \ref{stab}.}\\
Pick $A,B,\varphi$ as above and let 
$$C=\varphi(A), \qquad c_i=\varphi(a_i), \qquad G=\langle Iso(\U|A),Iso(\U|B)\rangle  \ .$$
Using the finite injectivity of $\U$ and the triangle inequality, we see that by moving each $c_i$ ($i>k$) a little bit we can assume that $c_i \not \in A\cup B$ for all $i>k$. \\
The idea of the proof is that, starting from any $\psi(a_1),\ldots,\psi(a_n)$ with $\psi\in \langle Iso(\U|A),Iso(\U|B)\rangle$ one can find some 
other $\tilde{\psi}(a_1),\ldots,\tilde{\psi}(a_n)$ with $\tilde{\psi}$ still in $\langle Iso(\U|A),Iso(\U|B)\rangle$ and $\tilde{\psi}(a_i)$ being closer to $c_i$ than $\psi(a_i)$.

To see why this is possible, one needs to perform various operations, which are detailed in a series of lemmas.
We first have to make a technical assumption whose usefulness will be apparent only later on; on first reading it may be a good idea to skip it and read the rest of the proof to see why it is useful.

\begin{lem} \label{depart}
Without loss of generality, one can assume that every triangle $\{a_p,c_q,b_r\}$ where $q>k$ and $\max(p,r)>k$ is flat. Similarly, one can assume that no triangle $\{b_p,c_q,c_r\}$ or $\{a_p,c_q,c_r\}$ with $p,q>k$ is flat. 
\end{lem}

\noindent{\bf Proof of Lemma \ref{depart}} \\
Let us introduce some notation. Denote by ${\mathcal C}$ the collection of metric extensions of $A \cup B$ by a set $\{x_1,\ldots,x_n\}$ isometric to $\{a_1,\ldots,a_n\}$ and such that $x_i=a_i$ for $i \le k $.\\ 
We may see $C$ as an element of this collection. Now, let $\Delta$ be one of the triangles this lemma is concerned with (say, $\Delta=\{a_p,c_q,b_r\}$). For $\overline{d}=\{d_1,\ldots,d_n\} \in {\mathcal C}$ we let $\Delta_{\overline{d}}$ denote the metric triangle obtained by replacing each $c_i$ in $\Delta$ by $d_i$ (so in our example $\Delta_{\overline{d}}=\{a_p,d_q,b_r\}$).\\
For each triangle $\Delta$, it is not too hard to see that there exists some $\overline{d} \in {\mathcal C}$ such that 
$\Delta_{\overline{d}}$ is not flat. So, taking the average of all these extensions, we see that there exists an element $\overline{e}$ of ${\mathcal C}$ such that each $\Delta_{\overline{e}}$ is not flat.\\
Finally, consider for each $\delta>0$ the element $\overline{e}^{\delta}$ of ${\mathcal C}$ obtained by setting 
$$\forall z \in A \cup B \ d(e^{\delta}_i,z)= (1-\delta)d(c_i,z)+\delta d(e_i,z) \ .  $$
This is an element of ${\mathcal C}$ such that none of the triangles we care about are flat; also, using Lemma \ref{lemlocal1} and the finite injectivity of $\U$, we see that we can assume that this extension is realized by points $\delta$-close to the $c_i$'s.\\
In other words, by moving each $c_i$ very slightly we can assume that each triangle we care about is not flat. Hence proving Proposition \ref{stab} in the case when all these triangles are not flat is enough to prove it altogether. \hfill $\square_{\text{\ref{depart}}}$\\

We now assume that all of the triangles mentioned above are not flat.
Denote by ${\mathcal O}$ the closure (in $\U^n$) of the set 
$$\{g(a_1),\ldots,g(a_n) \colon g \in G\}\ . $$
We also let $F(x_1,\ldots,x_n)= \sum d(x_i,c_i)$ for $(x_1,\ldots,x_n) \in {\mathcal O}$, and
$$R=\inf\{F(x_1,\ldots,x_n) \colon (x_1,\ldots,x_n) \in {\mathcal O}\} \ . $$ 
Note that ${\mathcal O}$ is $G$-invariant; what we want to prove is that $R=0$. We proceed by contradiction and assume that $R>0$. Let us begin by showing that $R$ is actually a minimum.

\begin{lem} \label{atteint}
There exists $(x_1,\ldots,x_n) \in {\mathcal O}$ such that $F(x_1,\ldots,x_n)=R$.
\end{lem}

\noindent{\bf Proof of Lemma \ref{atteint}.}\\
Pick a sequence $\overline{x^i}$ in $\mathcal{O}$ such that $F(\overline{x^i})$ converges to $R$. Note that up to some extraction one may assume that $d(x^i_j,a_p)$ and $d(x^i_j,c_p)$ converge, for all fixed $p$, when $i$ goes to infinity. \\
Hence, using Lemma \ref{lemlocal1} and the finite injectivity of $\U$, we see that there exist sequences $(y^i_j)$ for all $j\in \{1,\ldots,n\}$ such that
\begin{itemize}
 \item $d(y^i_j,a_p)=d(x^i_j,a_p)$ and $d(y^i_j,c_p)=d(x^i_j,c_p)$ for all $i$ and all $j,p \in \{1,\ldots,n\}$,
 \item $d(y^i_j,y^i_p)=d(x^i_p,x^i_p)$ ($=d(a_j,a_p)$)
 \item Each sequence $(y^i_j)$ converges to some $x_j$.
\end{itemize}
Now note that the conditions above mean in particular that for all $i$ one can map $\overline{x^i}$ to $\overline{y^i}$ by an isometry fixing $A$, hence each $\overline{y^i}$ belongs to ${\mathcal O}$. Since ${\mathcal O}$ is closed we obtain that $(x_1,\ldots,x_n) \in {\mathcal O}$, and then it is clear that $F(x_1,\ldots,x_n)=R$. \hfill $\square_{\text{\ref{atteint}}}$\\

Now fix $\overline{x} \in {\mathcal O}$ such that $F(\overline{x})=R$.

\begin{lem} \label{contra1}
Assume that there exists some $i_0$ such that 
$$\forall j \in \{1,\ldots,n\} \ |d(c_{i_0},b_j)-d(x_{i_0},b_j)|<d(c_{i_0},x_{i_0})$$
Then we reach a contradiction (and so it must be that $R=0$).
\end{lem}

\noindent{\bf Proof of Lemma \ref{contra1}.} \\
We know that if  $\{z_1,\ldots,z_n\}$ is an abstract extension of $B$ by a set isometric to $\{x_1,\ldots,x_n\}$ and such that $d(z_i,b_j)=d(x_i,b_j)$ for all $i,j$, then this extension is realized in $\U$ by points $z_1,\ldots,z_n$ satisfying additionally (thanks to lemma \ref{lemlocal1})
$$\forall i \ d(z_i,c_i)= \sup\{|d(c_i,b_j) - d(x_i,b_j)| \colon j=1,\ldots,n  \} \ .$$
Such $z_1,\ldots,z_n$ are equal to $g(x_1),\ldots,g(x_n)$ for some $g$ in $Iso(\U|B)$, 
which shows that $(z_1,\ldots,z_n) \in {\mathcal O}$.\\
We have by construction $d(z_i,c_i) \le d(x_i,c_i)$ for all $i$ and the only way that $d(z_{i_0},c_{i_0})$ is not strictly less than $d(x_{i_0},c_{i_0})$ is that there exists $j$ such that $|d(c_{i_0},b_j)-d(x_{i_0},b_j)|=d(c_{i_0},x_{i_0})$. This is ruled out by the hypothesis of the lemma, and so $F(z_1,\ldots,z_n)< F(x_1,\ldots,x_n)$, which is also impossible. So our assumption that $R>0$ must be wrong.$\hfill \square_{\text{\ref{contra1}}}$\\

We would have finished the proof if not for the possible existence of those pesky flat triangles, which we now must remove. That is the purpose of the following lemma.

\begin{lem} \label{contra2}
There exists $y_1,\ldots, y_n \in {\mathcal O}$ such that $F(y_1,\ldots,y_n)=R$ and for some $i_0$ one has $|d(c_{i_0},b_j)-d(y_{i_0},b_j)|<d(c_{i_0},y_{i_0})$ for all $j$.
\end{lem}

\noindent{\bf Proof of Lemma \ref{contra2}.} Pick $(x_1,\ldots,x_n) \in {\mathcal O}$ such that $F(x_1,\ldots,x_n)=R$ and fix $i_0>k$ such that $x_{i_0} \ne c_{i_0}$ (if there is no such $i_0$ then $R=0$, which is impossible by assumption).\\
Fix also $j_0$ such that $|d(c_{i_0},b_{j_0})-d(x_{i_0},b_{j_0})|=d(c_{i_0},x_{i_0})$. 

\begin{claim*}
There exist $z_1,\ldots,z_n \in \U$ such that 
\begin{itemize} 
\item[(1)] $d(z_i,z_j)=d(a_i,a_j)$ for any $i,j$;
 \item[(2)] $d(z_i,a_j)=d(x_i,a_j)$ for all $i$ and $j$;
\item[(3)] $d(z_i,c_i)=d(x_i,c_i)$;
\item[(4)] $|d(c_{i_0},b_{j_0})-d(z_{i_0},b_{j_0})|<d(c_{i_0},z_{i_0})$.
\end{itemize}
\end{claim*}
This claim is the heart of the lemma; indeed, it provides an abstract extension that removes one of the offending flat triangles, and then one just has to use a convex combination of all these extensions to obtain an abstract extension of 
$A \cup B \cup C$ by elements $y_1,\ldots,y_n$ which
satisfies conditions (1) through (3) above, and is such that no triangle $\{y_{i_0},c_{i_0},b_j\}$ is flat.\\
Using the finite injectivity of $\U$, one can assume that $y_1,\ldots,y_n$ are in $\U$, and then they must be in ${\mathcal O}$ since there is an isometry fixing $A$ and mapping $\{x_1,\ldots,x_n\}$ to $\{y_1,\ldots,y_n\}$ because of condition (2). Then condition (3) ensures that $F(x_1\ldots,x_n)=F(y_1,\ldots,y_n)$, which was the last point to check.\\

\noindent{\bf Proof of the claim.} \\
Assume first that $d(c_{i_0},b_{j_0})-d(x_{i_0},b_{j_0})=d(c_{i_0},x_{i_0})$.\\ 
Then we want to define an abstract extension of $A \cup B \cup C$ by a set isometric to $\{a_1,\ldots,a_n\}$ by setting
\begin{itemize} \item $d(z_{i_0},b_{j_0})= d(x_{i_0},b_{j_0})+\delta $
\item $d(z_{i_0},x) =  d(x_{i_0},x)$ if  $x \ne b_{j_0}$.
  \end{itemize}

In other words, we want to increase the distance to $b_{j_0}$ while keeping all others constant; if this is not possible even for very small 
$\delta >0$, then it must happen that for some $x\in A \cup B \cup C \setminus \{b_{j_0}\} $ one has 
$$d(x_{i_0},x)+d(x,b_{j_0})=d(x_{i_0},b_{j_0}) \ .$$
From this we get, using the triangle inequality (a picture is useful here!)
$$d(c_{i_0},x_{i_0})+d(x_{i_0},x)+d(x,b_{j_0})=d(c_{i_0},b_{j_0}) $$
This implies in turn that 
\begin{eqnarray}
d(c_{i_0},x_{i_0})+d(x_{i_0},x)\  = & d(c_{i_0},x)\\
d(c_{i_0},x)+d(x,b_{j_0})\  =& d(c_{i_0},b_{j_0})
\end{eqnarray}

If $x$ does not belong to $A\cap B$, or if $j_0 \ge k+1$, then (2) is prohibited by Lemma \ref{depart}.\\
So it must be that $x$ belongs to $A\cap B$ and $j_0\le k$. But then one must have 
$d(c_{i_0},x)=d(x_{i_0},x)$, and from (1) we then get $c_{i_0}=x_{i_0}$. This contradicts our choice of $i_0$.\\

The case when $d(x_{i_0},c_{i_0})+d(c_{i_0},b_{j_0})=d(x_{i_0},b_{j_0})$ is similar, except this time one wants to decrease the distance to $b_{j_0}$ while keeping all others constant. If it's not possible then the same line of reasoning as above leads to the desired contradiction. \hfill $\square_{\text{\ref{contra2}}}$\\

This concludes the proof of Proposition \ref{stab}. \hfill $\square_{\text{\ref{stab}}}$\\

\begin{thm} \label{Isaac}
Let $A,B \subseteq \U$ be finite sets. Then
$$Iso(\U |A\cap B)= \langle Iso(\U|A),Iso(\U|B)\rangle  $$
\end{thm}

\noindent{\bf Proof of Theorem \ref{Isaac}.}\\
Let $A,B \subseteq \U$ be finite and $\varphi \in Iso(\U |A\cap B)$. Let $V$ be a basic neighborhood of $\varphi$, which we can assume without loss of generality to be of the form 
$$ V=\{\psi \in Iso(\U) \colon d(\psi(x_i),\varphi(x_i)) < \varepsilon \} $$
with $\varepsilon>0$, $A \subseteq \{x_1,\ldots,x_n\}=X_0$.

\begin{claim*}
Without loss of generality, one can assume that $A\cap B=X_0 \cap B$. 
\end{claim*}
 
\noindent{\bf Proof.} We have to consider the case when some points of $X_0$ belong to $B \setminus A$; so we are concerned with points in $(X_0 \setminus A) \cap B$. Set $\varepsilon'=\varepsilon/3$. Since $(X_0 \setminus A) \cap B$ is finite, there exists for any $x_j \in (X_0 \setminus A) \cap B$ a point $x'_j$ which does not belong to $B$ and is such that $d(x_j,x'_j)< \varepsilon'$. Replace each $x_j$ by $x'_j$, and denote by $X'_0$ the set thus obtained. It still contains $A$; furthermore, the open set
$$V'=\{\psi \in Iso(\U) \colon \forall x \in X'_0\ d(\psi(x),\varphi(x)) < \varepsilon'\} $$
is contained in $V$. So, replacing $V$ by $V'$ and $X_0$ by $X'_0$ we are now in the situation described by the claim.  \hfill $\square_{\text{claim}}$\\

Using an amalgamation over $A \cap B$, one can increase $X_0$ to some finite set $X$ and $B$ to some finite set $\tilde X$ isometric to $X$ (by an isometry fixing $A \cap B$) such that additionally $X \cap \tilde X=A \cap B$.\\
Then we can find some isometry $\psi$ with $d(\psi(x_i),\varphi(x_i))< \varepsilon/2$ for all $i$, $\psi(y)=y$ for all $y \in X \cap \tilde X$, and $\psi(X \setminus \tilde X) \cap \tilde X= \emptyset$. \\
Applying Proposition \ref{stab} to $X,\tilde X,\psi$ we get an isometry 
$$\tilde{\psi} \in \, \langle Iso(\U|X),Iso(\U|\tilde X)\rangle  \,\subseteq \, \langle Iso(\U|A),Iso(\U|B)\rangle $$ 
such that $d(\tilde{\psi}(x_i),\psi(x_i)) < \varepsilon/2$, so the triangle inequality gives us that $\tilde{\psi}$ is in $V$. This proves that 
$\langle Iso(\U|A),Iso(\U|B)\rangle $ is dense in $Iso(\U |A\cap B)$, and this concludes the proof. $\hfill \square_{\text{\ref{Isaac}}}$\\

When seeing the statement of theorem \ref{Isaac}, it is natural to wonder whether one really needs to consider the closure of the subgroup $H$ generated by $Iso(\U|A)$ and $Iso(\U|B)$: it may be that one already has $H=Iso(\U|A\cap B)$ and no closure operation is needed. 
To see that it is needed indeed, consider the case when $A,B$ are nonempty finite subsets of $\U$ with empty intersection. In this case $Iso(\U|A \cap B)=Iso(\U)$.\\
For any isometry $\varphi$ that fixes either $A$ or $B$ one must have, because of the triangle inequality:
$$\forall x \in \U \  d(x,\varphi(x)) \le 2d(x,A\cup B )\ . $$
Now consider a product of $n$ isometries $\varphi_1,\ldots,\varphi_n$ belonging to 
$Iso(\U|A) \cup Iso(\U|B)$, and observe that
$$\forall x \in \U \ d(\varphi_1\ldots\varphi_n(x),x) \le d(\varphi_2\ldots\varphi_n(x),x)+d(\varphi_1(x),x)\ .$$
By induction, we obtain that any isometry $\psi$ which can be written as a product of $n$ elements of 
$Iso(\U|A) \cup Iso(\U|B)$ must satisfy
$$\forall x \in \U  \ d(x,\varphi(x)) \le 2 n d(x,A\cup B)\ . $$
The set of isometries satisfying those conditions is meager in $Iso(\U)$, and since this is true for all $n$, $H$ itself is meager in $Iso(\U)$.\\

We conclude this section, and this article, by pointing out that the ideas presented here may be used to prove that Theorem \ref{Isaac} is still true when one replaces $\U$ by the Urysohn space of diameter $1$ (which is the unique, up to isometry, Polish metric space of diameter $1$ which is both ultrahomogeneous and universal for separable metric spaces of diameter $1$). Note also that it seems reasonable to expect that the same theorem holds with $A,B$ compact instead of finite, but I did not try to check the details.\\

\noindent \emph{Acknolewdgment.} I'm very grateful to the anonymous referee for a careful reading of the first version of this paper and several valuable observations and suggestions.

\begin{small}
\noindent Julien Melleray \\
melleray@math.univ-lyon1.fr \\

\noindent Universit\'e de Lyon; \\
Universit\'e Lyon 1; \\
INSA de Lyon, F-69621; \\
Ecole Centrale de Lyon; \\
CNRS, UMR 5208, Institut Camille Jordan,\\
43 blvd du 11 novembre 1918, \\
F-69622 Villeurbanne Cedex, France.
\end{small}

\end{document}